\documentclass{article}
\usepackage{amsmath}
\newtheorem{theorem}{Theorem}[section]
\newtheorem{proposition}{Proposition}[subsection]
\newtheorem{lemma}[theorem]{Lemma}

\newtheorem{remark}[theorem]{Remark}

\newenvironment{proof}{{\bf Proof.}}{\par\hspace{25em}\rule{1ex}{1ex}\par}
\title{On the mean curvature of a unit vector
field.\footnote{Publ. Math. Debrecen 2002 60, 2/3, 131-155}}
\author{Yampolsky A.\\[1ex]
Department of Mechanics and Mathematics, \\
    Kharkov National University, Ukraine.\\
    e-mail: ALYampolsky@univer.kharkov.ua}
\date{}
\begin{document}
\maketitle
\begin{abstract}
    We present an explicit formula for the mean curvature of a unit
    vector field on a Riemannian manifold, using a special but natural
    frame. As applications, we treat some known and new examples of
    minimal unit vector fields. We also give an example of a vector
    field of constant mean curvature on the Lobachevsky $(n+1)$ space.\\[1ex]
    {\it AMS subject class:} Primary 54C40,14E20; Secondary 46E25, 20C20
\end{abstract}

\section*{Introduction}

    Let $(M,g)$ be an $n+1$ -- dimensional Riemannian manifold  with
    metric $g$. A vector field $\xi$ on it is called {\it holonomic}
    if $\xi$ is a field of normals of some family of regular
    hypersurfaces in $M$ and {\it non-holonomic} otherwise. The
    foundation of the classical geometry of unit vector fields was
    proposed by A.Voss at the end of the nineteenth century. The
    theory includes the {\it Gaussian} and {\it the mean curvature} of
    a vector field and their generalizations (see \cite{Am} for
    details). Here we will consider a unit vector field from another
    point of view. Namely, let $T_1M$ be the unit tangent sphere
    bundle of $M$ endowed with the Sasaki metric \cite{S}. If $\xi$ is
    a unit vector field on $M$, then one may consider $\xi$ as a
    mapping $\xi : M \to T_1M $ so that the image $\xi (M) $  is a
    submanifold in $T_1M$ with the metric induced  from $T_1M$.
    H.Gluck and W.Ziller \cite{G-Z} called $\xi$ {\it a minimal vector
    field} if $ \xi(M)$ is  of minimal volume with respect to induced
    metric. They considered the unit vector field  on $S^3$ tangent to
    the fibers of a Hopf fibration $S^3
    \stackrel{S^1}{\longrightarrow} S^2$ and proved that these (Hopf)
    vector fields are unique ones with global minimal volume. Note
    that this result is not true for greater dimensions where Hopf
    vector fields are still critical points for the volume functional
    but do not provide the global minimum among all unit vector fields
    \cite{Jon,Ped}. The local aspect of the problem was considered
    first in \cite{GM-LF}. The authors have found the necessary and
    sufficient condition for a unit vector field to generate locally a
    minimal submanifold in the tangent sphere bundle. In fact,  that
    condition implies that {\it the mean curvature } of the
    submanifold $\xi(M)$  is zero. Using that criterion, a number of
    examples of local minimal vector unit fields have been found (~see
lab2    \cite{ BX-V1, BX-V2, GD-V1, GD-V2, TS-V1,TS-V2}).

    In this paper, we give an {\it explicit formula} for the mean curvature of
    $\xi(M)$ using some special but natural normal frame for $\xi(M)$ and give an
    example of a unit vector field of {\it constant mean curvature}
    on a Lobachevsky space. We shall state the main result after some preliminaries.

    Let $\nabla$ denote the Levi-Civita connection on $M$. Then $\nabla_ X \xi$ is
    always orthogonal to $\xi$ and hence,
    $(\nabla\xi)(X)=\nabla_X\xi :T_pM \to \xi^\perp_p$ is a linear operator at each
    $p\in M$. We define the adjoint operator $ (\nabla\xi)^*(X) :\xi^\perp_p \to
    T_pM$  by
    $$
        \left< (\nabla\xi)^*X,Y\right>_g = \left< X,\nabla_Y\xi \right>_g
    $$
    Then there is an orthonormal frame $e_0, e_1, \dots , e_n $ in $T_pM$ and an
    orthonormal frame $f_1, \dots , f_n $ in $\xi_p^\perp$ such that
    $$
        (\nabla\xi)(e_0)=0 , \quad (\nabla\xi)(e_\alpha )=\lambda_\alpha f_\alpha,
        \quad (\nabla\xi)^*(f_\alpha)=\lambda_\alpha e_\alpha ,
        \qquad \alpha=1, \dots , n ,
    $$
    where $\lambda_1\ge \lambda_{2}\ge \dots \ge \lambda_n\ge 0$ are the singular
    values of $\nabla\xi$. As we will see, the vectors
    $$
    \tilde n_{\sigma |} =\frac{1}{\sqrt{1+\lambda_\sigma^2}}\big(-\lambda_\sigma
    e_\sigma^h +f_\sigma^v \ \big),\mbox{\hspace{3em}} \sigma=1,\dots , n ,
    $$
    where $H$ and $V$ are the horizontal and vertical lifts
    respectively, form  an orthonormal frame in the normal bundle of $\xi(M)$.

    Furthermore, we introduce the notation
    $$
        r(X,Y)\xi=\nabla_X\nabla_Y\xi-\nabla_{\nabla_XY}\xi.
    $$
    Then $R(X,Y)\xi=r(X,Y)\xi-r(Y,X)\xi$ , where $R$ is the Riemannian curvature
    tensor. Now we are able to state our main result.
\vspace{1ex}

    {\bf Theorem \ref{Th1}} {\it Let $H_{\sigma |}$ be the components of
    the mean curvature vector of $\xi(M)$ with respect to the
    orthonormal frame $\tilde n_\sigma$. Then
$$
\begin{array}{c}
    (n+1)H_{\sigma |}=\\[1ex]\displaystyle
    \frac{1}{\sqrt{1+\lambda_\sigma^2}}
    \left\{ \big<r(e_0,e_0)\xi,f_\sigma\big> + \sum_{\alpha =1}^n
    \frac{\big<r(e_\alpha,e_\alpha)\xi,f_\sigma\big>+
    \lambda_\sigma\lambda_\alpha \big<R(e_\sigma,
    e_\alpha)\xi,f_\alpha)\big>}{1+\lambda_\alpha^2} \right\}.
\end{array}
$$ }. \vspace{1ex}

    The following very simple example gives a unit vector field of
    {\it constant mean curvature}.
\vspace{1ex}

    {\bf Proposition \ref{Ex}} {\it Let $M$ be the Lobachevsky 2-plane
    with the metric
    $$
     ds^2=du^2+e^{2u}dv^2.
    $$
    Let $X_1=\{1,0\}$ and $X_2=\{0,e^{-u}\}$. Then $\xi=\cos \omega X_1+\sin\omega X_2$,
    where $\omega=au+b $, generates a hypersurface $\xi(M)\subset T_1M$ of constant
    mean curvature
    $$
    H=\frac{a}{2\sqrt{2+a^2}}.
    $$
    }
\vspace{1ex}

    {\bf Index convention.} Throughout the paper we take $i, j, k, \ldots =0, \dots, n$
    and $\alpha, \beta, \ldots = 1, \dots, n .$

%\newpage
\section{ Basic concepts from the geometry of the unit tangent sphere bundle.}

    Let $(u^0,\dots ,u^n)$ be a local coordinate system on $M$ and
    let $\partial /\partial u^i $ be the vectors of a natural frame on
    $M^n.$ The points of the tangent bundle  $TM$ are the pairs
    $\tilde Q=(Q,\xi)$, where $Q\in M$ and $\xi\in T_QM$.  Each point
    $\tilde Q\in TM$ is uniquely determined by the set of parameters
    $(u^0,\dots ,u^n;\xi^0,\dots ,\xi^n)$, where $(u^0, \dots ,u^n)$
    fix the point $Q$ and  $\{\xi^0, \dots, \xi^n\}$  are the
    coordinates of $\xi$ with respect to the frame $\{ \partial
    /\partial u^0, \dots ,\partial /\partial u^n \}$. The local
    coordinates  $(u^0,\dots ,u^n;\xi^0, \dots ,\xi^n)$ are called
    {\it natural induced coordinates } in the tangent bundle. Each
    smooth tangent vector field $\xi=\xi(u^0, \dots ,u^n)$ generates a
    smooth submanifold $\xi(M)\subset TM$ having a parametric
    representation of the form
\begin{equation}
\label{lab1}
    \left\{
        \begin{array}{lcl}
                u^i    &  =  & u^i,\\
                \xi^i  &  =  & \xi ^i(u^0, \dots, u^n).
        \end{array}
    \right.
\end{equation}
    Setting $|\xi |=1$, we get a submanifold in the unit tangent sphere bundle
    $\xi(M^n)\subset T_1M^n.$

    A natural Riemannian metric on the tangent bundle has been defined by S.Sasaki
    \cite{S}. We describe it in terms of the {\it connection map}.

    The tangent space $T_{\tilde Q}TM$ can be split into {\it vertical} and {\it
    horizontal} parts:
    $$
    T_{\tilde Q}TM^n=H_{\tilde Q}TM^n \oplus V_{\tilde Q}TM^n.
    $$
    The vertical part  $V_{\tilde Q}TM$ is tangent to the fiber, while the
    horizontal part is transversal to it. For $\tilde X \in T_{\tilde Q}TM^n$ we
    have
\begin{equation}
    \label{lab2}
        \tilde X=\tilde X^i \partial /\partial u^i + \tilde X^{n+i} \partial /\partial
        \xi^i
\end{equation}
    with respect to the natural frame $\{ \partial /\partial u^i, \partial /\partial
    \xi^i \}$ on $TM$.

    Let  $\pi:TM \to M$ be the  projection map. It is easy to check that the
    differential $\pi_*:T_{\tilde Q}TM \to T_QM $ of the mapping $\pi$  acts on
    $\tilde X$ as follows:
\begin{equation}
    \label{lab3}
        \pi_*\tilde X=\tilde X^i \partial /\partial u^i,
\end{equation}
    and is a linear isomorphism between $V_{\tilde Q}TM$ and $T_QM$.

    The {\it connection map} $K: T_{\tilde Q}TM \to T_QM$ acts on $\tilde X$ by
\begin{equation}
    \label{lab4}
        K\tilde X=(\tilde X^{n+i}+\Gamma_{jk}^{i}\xi^j\tilde X^k) \partial /\partial u^i
\end{equation}
    and it is a linear isomorphism between $H_{\tilde Q}TM$ and $T_QM$.
    Moreover, it is easy to see that $V_{\tilde Q}TM=\ker \pi_*$,  $H_{\tilde Q}TM=\ker K$.
    The  images $\pi_*\tilde X$ and  $K\tilde X$ are called {\it horizontal} and
    {\it vertical } projections of  $\tilde X$, respectively.

    The {\it Sasaki metric} on $TM$ is defined by the following scalar product:
    if $\tilde X,\tilde Y \in T_{\tilde Q}TM$,  then
\begin{equation}
\label{lab5}
        \big<\big< \tilde X,\tilde Y \big>\big>_S=
        \big<\pi_* \tilde X, \pi_* \tilde Y\big>_g+\big<K \tilde X,K \tilde Y\big>_g
\end{equation}
    where  $\big<,\big>_g$ is the scalar  product with respect to the metric  $g $
    on the initial manifold  (the base space of tangent bundle).
    Horizontal and vertical subspaces are mutually orthogonal with respect to Sasaki
    metric.

    The inverse operations of projections (\ref{lab3}) and (\ref{lab4})
    are called {\it lifts}. Namely, if $X \in T_QM^n$, then
    $$
    X^H=X^i \partial /\partial u^i -\Gamma_{jk}^i\xi^j X^k \partial /\partial \xi^i
    $$
    is in $H_{\tilde Q}TM$ and is called the {\it horizontal lift } of  X, and
    $$
    X^V=X^i \partial /\partial \xi^i
    $$
    is in $V_{\tilde Q}TM$) and is called the {\it vertical lift } of  $ X$.

    Among all lifts of various vectors from  $T_QM$ into $T_{(Q,\xi)}TM$, one can
    naturally distinguish two of them, namely $\xi^H$ and $\xi^V$. The vector field
    $\xi^H$ is the {\it geodesic flow} vector field, while $\xi^V$ (being
    normalized) is a {\it unit normal} vector field of  $T_1M \subset TM$.

    In the geometry of the {\it unit tangent sphere bundle} it appears to be
    convenient to introduce the notion of {\it tangential lift} \cite{BX-V3}:
\begin{equation}\label{lab5_1}
    X^t=X^V-\big<X,\xi\big>\xi^V.
\end{equation}
    In other words, the tangential lift is the projection of the vertical lift onto the
    tangent space of $T_1M$.

    We denote by $\tilde\nabla$ the Levi-Civita connection of the Sasaki metric on
    $T_1M$. In terms of horizontal and tangential lifts we then have \cite{BX-V3}:
\begin{equation}
    \label{lab6}
    \begin{array}{ll}
        \tilde\nabla_{X^H}Y^H  =  (\nabla_XY)^H - \frac{1}{2}(R(X,Y)\xi)^t,
         &\tilde\nabla_{X^t}Y^H  =  \frac{1}{2}(R(\xi,X)Y)^H, \\
        \tilde\nabla_{X^H}Y^t  =  (\nabla_XY)^t \ + \frac{1}{2}(R(\xi_1,Y)X)^H, &
        \tilde\nabla_{X^t}Y^t  =  -\big<Y,\xi\big>X^t.
    \end{array}
\end{equation}
\begin{remark} \rm It is evident that if $Z \perp  \xi $,  the vertical and tangential
    lifts of  $Z$ coincide, particulary $ (\nabla_X\xi)^t=(\nabla_X\xi)^V$ for any
    $X$.  We will use this fact throughout the paper without special comments.
\end{remark}

\section{The mean curvature formula for a unit vector field}
\subsection{The structure of tangent and normal bundles of $\xi(M)$}

    Let $\xi$  be the unit tangent vector field on $M$. We denote by $T\xi(M)$ the
    tangent bundle of $\xi(M)\subset T_1M$. The structure of $T\xi(M)$ can be
    described as follows:

\begin{lemma}
    \label{L1} \it The vector  $\tilde X \in T_{(Q,\xi)}T_1M$ is tangent to $\xi(M)$
    at $(Q,\xi)$ if and only if
    \begin{equation}
        \label{lab7}
        \tilde X = X^H + (\nabla_X\xi)^V
    \end{equation}
    where $X \in T_QM$.
\end{lemma}

\begin{proof}  Using the local representation (\ref{lab1}) of $ \xi (M)$, we
    consider the coordinate frame of $T_{(Q,\xi)}\xi (M)$:
    $$
    \tilde e_i = \left\{0,\dots , 1, 0,  \dots , 0; \frac{\partial \xi^0}{\partial
    u^i},
    \dots , \frac{\partial \xi^n}{\partial u^i}\right\}.
    $$
    Let  $\tilde X \in T_{(Q,\xi)}TM$ be tangent to $\xi(M)$.
    Then
    $$
    \tilde X = \tilde X^i \tilde e_i.
    $$
    Applying  (\ref{lab3}) and  (\ref{lab4}), we obtain
    $$
    \begin{array}{lcl}
    \pi_*\tilde e_i  & = & \partial /\partial u^i,  \\
    K\tilde e_i      & = & \nabla_i\xi.
    \end{array}
    $$
    From this we get
    $$
    \begin{array}{lcl}
    \pi_* \tilde X & = & \tilde X^i \partial /\partial u^i, \\
    K \tilde X     & = & \nabla_{\pi_* \tilde X} \xi.
    \end{array}
    $$
    Setting  $X = \pi_* \tilde X$ and taking into account the remark, we get
    (\ref{lab7}).

\end{proof}

    To describe the structure of the normal bundle of  $\xi(M)$, we use the {\it
    adjoint covariant derivative operator}. As $\xi $ is a fixed unit vector field,
    $\nabla_X\xi$  can be considered as a pointwise linear operator
    $(\nabla\xi):T_QM \to \xi^{\perp}$, where $\xi^{\perp}$  is the orthogonal
    complement  of $\xi$ in $T_QM$, acting as
    $$
     (\nabla \xi)(X) = \nabla _X \xi.
    $$
    The matrix of this operator is formed by the covariant derivatives $\nabla_i
    \xi^k$.

    The {\it adjoint covariant derivative} linear operator $(\nabla\xi)^*:
    \xi^{\perp} \to T_QM$ can be defined in a standard way:
\begin{equation}
    \label{lab8}
        \big<(\nabla \xi)^*X,Y\big> = \big<X,(\nabla\xi)(Y)\big>
\end{equation}
    for each $X \in \xi^{\perp}$. The matrix of $(\nabla \xi)^*$ has the form
    $$
    \left[ (\nabla \xi)^* \right]_j^i = g^{im} \nabla _m \xi^kg_{kj}.
    $$

    As $\nabla$ is the Riemannian connection for $g$,  we obtain for $(\nabla\xi)^*$
    the formally transposed matrix
    $$
    \left[(\nabla \xi)^* \right]_k^i = \nabla^i \xi_k.
    $$

    Now the structure of $\xi(M)$ can be described as follows:

\begin{lemma}
    \label{L2}
         The vector $\tilde N \in T_{(Q,\xi)}T_1M$ is  normal to
        $\xi(M)$ if and only if $$ \tilde N = - \left[(\nabla
        \xi)^*N\right]^H + N^V $$ where $N \in T_QM$ and $N\perp \xi$.
\end{lemma}

    The proof follows easily from (\ref{lab5}), (\ref{lab7}) and (\ref{lab8})

\subsection{Second fundamental form of $\xi(M)$ in $T_1M$}

    We denote by $\tilde \Omega_{\tilde N}$  the second fundamental form of $\xi(M)$
    in $T_1M^n$ with respect to the normal vector field $\tilde N$ defined in Lemma
    \ref{L2}.
    Then the following statement holds.

\begin{lemma} \label{L3}
     For  $\tilde X, \tilde Y $ being tangent to $\xi(M)$ we have
     $$
    \tilde \Omega_{\tilde N}(\tilde X, \tilde Y) = \frac{1}{2}
    \big< r(X,Y) \xi + r(Y,X) \xi -\nabla_{R(\xi, \nabla_X
    \xi)Y+R(\xi, \nabla_Y \xi)X} \xi ,N \big>,
    $$
      where $r(X,Y)\xi =\nabla_X \nabla_Y \xi - \nabla_{\nabla_X Y}\xi$
\end{lemma}

\begin{proof}
    By definition we have
     $$
     \tilde \Omega_{\tilde N}(\tilde X,\tilde
    Y) = \big<\big< \tilde \nabla_{\tilde X}\tilde Y,\tilde
    N\big>\big>
    $$
    where $\tilde X,\tilde Y \in T_{(Q,\xi)} \xi(M)$.
    Using Lemma \ref{L1}, we put $\tilde X = X^H + (\nabla_X\xi)^V;  \
    \tilde Y = Y^H +(\nabla_Y\xi)^V$. Then applying (\ref{lab6}) and
    (\ref{lab5_1}), we
    have
    $$
    \begin{array}{l}
    \tilde \nabla_{\tilde X} \tilde Y  =  \tilde
    \nabla_{X^H+(\nabla_X\xi)^t} (Y^H +(\nabla_Y\xi)^t) = \\[1ex]
    \left[\nabla_X Y + \frac{1}{2}R(\xi,\nabla_X \xi)Y +
    \frac{1}{2}R(\xi,\nabla_Y \xi)X\right]^H+
    \left[\nabla_X \nabla_Y \xi - \frac{1}{2}R(X,Y) \xi\right]^t= \\[1ex]
    \left[\nabla_X Y + \frac{1}{2}R(\xi,\nabla_X \xi)Y +
    \frac{1}{2}R(\xi,\nabla_Y \xi)X\right]^H + \left[\nabla_X \nabla_Y \xi -
    \frac{1}{2}R(X,Y) \xi\right]^V - \\[1ex]
    \big<\nabla_X \nabla_Y \xi,\xi\big>\xi^V.
    \end{array}
    $$
    Let $N$ be orthogonal to $\xi$. Then $\tilde N = - \left[
    (\nabla \xi)^*N \right]^H + N^V $ is normal to $\xi(M)$.
     Therefore

    \begin{eqnarray}
    \tilde \Omega_{\tilde N}(\tilde X, \tilde Y)  = -\big<\nabla_XY +
    \frac{1}{2}R(\xi,\nabla_X\xi)Y + \frac{1}{2}R(\xi,\nabla_Y\xi)X,
    (\nabla \xi)^*N\big> +\nonumber\\[1ex]
    \big<\nabla_X \nabla_Y \xi - \frac{1}{2}R(X,Y)\xi,N\big>=\nonumber\\[1ex]
    \big<\nabla_X \nabla_Y \xi - \frac{1}{2}R(X,Y)\xi
    -\nabla_{\nabla_XY+\frac{1}{2}R(\xi,\nabla_X \xi)Y+
    \frac{1}{2}R(\xi,\nabla_Y \xi)X} \xi, N\big>.\hspace{1em}\label{lab9}
   \end{eqnarray}

     To simplify the expression (\ref{lab9}), we introduce the
    following tensor $r$:
\begin{equation}\label{lab10}
    r(X,Y)\xi = \nabla_X \nabla_Y \xi - \nabla_{\nabla_X Y}\xi.
\end{equation}
    Then for the Riemannian tensor, we get
    $$
    R(X,Y)\xi = r(X,Y)\xi - r(Y,X)\xi
    $$
    and (\ref{lab9}) can be rewritten as
\begin{equation}\label{lab11}
    \tilde \Omega_{\tilde N}(\tilde X, \tilde Y) = \frac{1}{2} \big<
    r(X,Y) \xi + r(Y,X) \xi -\nabla_{R(\xi, \nabla_X \xi)Y+R(\xi,
    \nabla_Y \xi)X} \xi ,N \big>.
\end{equation}
\end{proof}

    Next, we determine the components of $\tilde \Omega$ with respect to some special
    frame.

    As $(\nabla \xi): T_QM \to \xi^\perp$ and $(\nabla \xi)^*:
    \xi^\perp \to T_QM$ are mutually adjoint, then in
    $T_QM$ and  $\xi^\perp$, respectively, there exist orthonormal frames
    $\{e_0, e_1,\dots ,e_n\}$ and $\{f_1, \dots ,f_n\}$ such that
    $$
    \left\{
    \begin{array}{lrl}
        (\nabla \xi)e_0 & = & 0, \\
        (\nabla \xi)e_\alpha & = & \lambda_\alpha f_\alpha, \\
        (\nabla \xi)^*f_\alpha & = & \lambda_\alpha e_\alpha,
    \end{array}
    \right.
    $$
    where $\lambda_n \ge \lambda_{n-1} \dots \ge \lambda_1 \ge 0$ is a set of
    singular values (functions) of the linear operator $\nabla \xi$.
    Then
\begin{equation}
\label{lab9'}
    \left\{
    \begin{array}{l}
\tilde e_0  =  e_0^H, \\ \tilde e_\alpha =
e_\alpha^H+(\nabla_{e_\alpha}\xi)^V = e_\alpha^H + \lambda_\alpha
f_\alpha^V
     \end{array}   \right.
\end{equation}
    form an orthogonal frame of the tangent space of
    $T_{(Q,\xi)}\xi(M)$ while
\begin{equation}\label{n12}
    \tilde n_{\sigma} = \frac{1}{\sqrt{1+\lambda_
    \sigma^2}}\left(\lambda_\sigma e_\sigma^H - f_\sigma^V \right)
\end{equation}
    form the orthonormal frame in  $\xi(M)^\perp$.

\begin{lemma}\label{L4}
    The components of second fundamental form of
    $\xi(M)\subset T_1M$ with respect to the frames (\ref{lab9'}) and
    (\ref{n12}) are given by
    $$
\begin{array}{rcl}
    \tilde \Omega_{\sigma | 00} &=&
    \frac{1}{\sqrt{1+\lambda_\sigma^2}} \big\{ \big< r(e_0,e_0)
    \xi,f_\sigma \big> \big\}, \\[1ex]
    \tilde \Omega_{\sigma | \alpha
    0} &=& \frac{1}{2}\frac{1}{\sqrt{1+\lambda_\sigma^2}} \frac{1}{\sqrt{1+\lambda_\alpha^2}}
    \big\{ \big<
    r(e_\alpha,e_0) \xi + r(e_0,e_\alpha) \xi,f_\sigma \big> +
    \lambda_\sigma \lambda_\alpha \big< R(e_\sigma,e_0) \xi, f_\alpha
    \big> \big\},   \\[1ex]
    \tilde \Omega_{\sigma | \alpha \beta}
    &=& \frac{1}{2}\frac{1}{ \sqrt{1+ \lambda_\sigma^2}} \frac{1}{\sqrt{1+\lambda_\alpha^2}}
    \frac{1}{\sqrt{1+\lambda_\beta^2}}
    \big\{ \big<
    r(e_\alpha, e_\beta) \xi+ r(e_\beta, e_\alpha) \xi, f_\sigma \big>\\
    &&+ \lambda_\alpha \lambda_\sigma \big< R(e_\sigma, e_\beta)
    \xi, f_\alpha \big> +  \lambda_\beta \lambda_\sigma \big<
    R(e_\sigma, e_\alpha) \xi, f_\beta \big> \big\},
\end{array}
    $$
    where $\sigma,\alpha,\beta=1,\dots,n$
\end{lemma}

\begin{proof}

    Indeed, with respect to (\ref{lab9'}) and (\ref{n12}) the
    components of $\tilde\Omega $ are $$ \tilde \Omega_{\sigma |
    ik} = \tilde \Omega_{\tilde n_\sigma}(\tilde e_i, \tilde e_k). $$

    Using (\ref{lab11}), we have

    $$ \tilde \Omega_{\sigma | ik}
    =\frac{1}{2}\frac{1}{\sqrt{1+\lambda_\sigma^2}} \big<r(e_i,e_k)
    \xi+r(e_k,e_i)\xi -\nabla_{R(\xi, \nabla_{e_i} \xi)e_k+R(\xi,
    \nabla_{e_k} \xi)e_i} \xi ,f_{\sigma} \big>. $$
     Setting $i=k=0$ and applying (\ref{lab9'}), we get
     $$
     \tilde \Omega_{\sigma | 00}= \frac{1}{\sqrt{1+\lambda_\sigma^2}}\big\{
      \big< r(e_0,e_0)\xi,f_\sigma \big> \big\}.
    $$
    Setting $i= \alpha,\, k = 0$ and applying (\ref{lab9'}) again, we obtain
    $$
\begin{array}{rl}
    \tilde \Omega_{\sigma | \alpha 0} = & \frac{1}{2}
    \frac{1}{\sqrt{1+ \lambda_\sigma^2}} \big\{ \big< r(e_\alpha, e_0)
    \xi, f_\sigma \big> + \big< r(e_0, e_\alpha) \xi, f_\sigma \big> -
    \big< \nabla_{R(\xi,( \nabla \xi)e_\alpha) e_0} \xi, f_\sigma
    \big> \big \} =\\[1ex]
    & \frac{1}{2} \frac{1}{\sqrt{1+\lambda_\sigma^2}} \big\{ \big<
    r(e_\alpha,e_0) \xi, f_\sigma \big> + \big< r(e_0,e_\alpha)
    \xi,f_\sigma \big> + \lambda_\sigma \lambda_\alpha \big<
    R(e_\sigma,e_0) \xi, f_\alpha \big> \big\}.
\end{array}
    $$
    Finally, setting $i=\alpha ,\, k=\beta $ applying  again (\ref{lab9'}), we obtain
    $$
    \begin{array}{lrl}
        \tilde \Omega_{\sigma | \alpha \beta} =
        &\frac{1}{2}\frac{1}{\sqrt{1+\lambda_\sigma^2}} &\left\{\big<
        r(e_\alpha, e_\beta) \xi + r(e_ \beta, e_\alpha) \xi - \right.\\
        & &\left.\quad\quad\quad\quad\quad \nabla_{R(\xi,( \nabla \xi)(e_\alpha))e_\beta +R(\xi,
        (\nabla \xi)(e_ \beta))e_\alpha} \xi, f_\sigma \big> \right\}=
        \\[1ex] & \frac12\frac{1}{\sqrt{1+\lambda_\sigma^2}} &\left\{\big<
        r(e_\alpha, e_\beta) \xi + r(e_\beta, e_\alpha) \xi,f_\sigma \big>
        - \right.\\
        &&\left.\quad\quad\big< \lambda_\alpha R(\xi,
        f_\alpha)e_\beta + \lambda_\beta R(\xi, f_\beta) e_\alpha,(\nabla
        \xi)^*(f_\sigma)\big> \right\}= \\[1ex]
        &\frac12\frac{1}{\sqrt{1+\lambda_\sigma^2}} &\left\{\big<
        r(e_\alpha,e_\beta) \xi + r(e_\beta, e_\alpha) \xi, f_\sigma \big>
        - \right.\\ &&\left.\quad\lambda_\alpha\lambda_\sigma \big<R(\xi,
        f_\alpha)e_\beta, e_\sigma \big> - \lambda_\beta\lambda_\sigma
        \big<R(\xi, f_\beta) e_\alpha, e_\sigma \big> \right\} = \\[1ex]
        &\frac{1}{2} \frac{1}{\sqrt{1+\lambda_\sigma^2}} &\left\{ \big<
        r(e_\alpha,e_\beta) \xi, f_\sigma \big> +\big< r(e_\beta,
        e_\alpha) \xi, f_\sigma \big> +\right.\\ &&\left.\quad
        \lambda_\alpha \lambda_\sigma \big< R(e_\sigma,e_\beta) \xi, f_
        \alpha \big> + \lambda_\beta \lambda_\sigma \big< R(e_\sigma,
        e_\alpha) \xi, f_\beta \big> \right\}.
    \end{array}
    $$

    So, the lemma is proved.
\end{proof}

\subsection{The mean curvature formula}

    Now we are able to prove the main result.
    \begin{theorem}\label{Th1}
    The components of the mean curvature vector of
    $\xi(M)\subset T_1M$ with respect to the frames (\ref{lab9'}) and
    (\ref{n12}) are given by
    \begin{equation}\label{H}
        \begin{array}{cc}
            (n+1)H_{\sigma |}=\\[1ex]\displaystyle
            \frac{1}{\sqrt{1+\lambda_\sigma^2}}\left\{
            \big<r(e_0,e_0)\xi,f_\sigma\big> + \sum\limits_{\alpha =1}^n
            \frac{\big<r(e_\alpha,e_\alpha)\xi,f_\sigma\big>+
            \lambda_\sigma\lambda_\alpha \big<R(e_\sigma,e_\alpha)\xi,
            f_\alpha)\big>}{1+\lambda_\alpha^2} \right\}.
        \end{array}
    \end{equation}
    \end{theorem}

\begin{proof} With respect to the frames (\ref{lab9'}) and (\ref{n12}) the
matrix of the first fundamental form  $\tilde G$ of $\xi(M)$ is

\begin{equation}
\label{n10} \tilde G = \left(
\begin{array}{cccc}
1             & 0              & \ldots      & 0       \\ 0
& 1+\lambda_1^2  & \ldots      &0   \\ \vdots        & \vdots
& \ddots      & \vdots  \\ 0             & 0              & \ldots
& 1+\lambda_{n}^2  \\
\end{array}
\right).
\end{equation}
For the inverse matrix we have

\begin{equation}
\label{n11} \tilde G^{-1} = \left(
\begin{array}{cccc}
1       & 0                      & \ldots           & 0  \\ 0
&\frac{1}{1+\lambda_1^2} &  \ldots          & 0  \\ \vdots  &
\vdots                 &\ddots            & \vdots \\ 0       & 0
&\ldots             & \frac{1}{1+\lambda_n^2}
\\
\end{array}
\right).
\end{equation}
So we have $$
\begin{array}{lcl}
\tilde \Omega_{\sigma | 00} & = &
 \frac{1}{\sqrt{1+\lambda_\sigma^2}} \big< r(e_0,e_0) \xi,f_\sigma \big>
,\\[1ex] \tilde \Omega_{\sigma | \alpha \alpha} & = &
 \frac{1}{\sqrt{1+\lambda_\sigma^2}} \big[ \big< r(e_\alpha,e_\alpha) \xi,
f_\sigma \big> + \lambda_ \sigma \lambda_ \alpha \big<
R(e_\sigma,e_\alpha) \xi, f_\alpha \big> \big].
\end{array}
$$
 Taking (\ref{n11})into account,  we have:
 $$ H_\sigma |
 =\frac{1}{(n+1)}\tilde G^{ii} \tilde \Omega_{\sigma | ii} = $$
$$ \frac{1}{(n+1)\sqrt{1+\lambda_\sigma^2}} \left\{ \big<
r(e_0,e_0)\xi,f_\sigma \big>+ \sum_{\alpha=1}^{n}  {\frac{ \big<
r(e_\alpha,e_\alpha) \xi,f_\sigma +\lambda_\sigma \lambda_\alpha
R(e_\sigma,e_\alpha) \xi,f_\alpha \big>} {1+\lambda_\alpha^2}}
\right\}. $$
 So we get the result.
\end{proof}
\subsubsection{Simplified formula for the mean curvature of a unit
    vector field.}
    It is possible to simplify the formula (\ref{H}). To do this, we introduce
    the following notations:
    $$
    E_{i| j k}=\big<\nabla_{\displaystyle e_i}e_j,e_k\big>, \quad
    F_{i| j k}=\big<\nabla_{\displaystyle e_i}f_j,f_k\big>,
    $$
    where $f_0$ is supposed to be zero. Evidently, $E_{i| j k}=-E_{i| kj}$
    and $F_{i| j k}=-F_{i| kj}$.
    Then it is simple to check that
    $$
    \big<r(e_i,e_j)\xi,f_k\big>=e_i(\lambda_j)\delta_{jk}+
    \lambda_jF_{i| jk}-\lambda_kE_{i| jk}.
    $$
    Therefore,
    $$
    \begin{array}{l}
    \big<r(e_j,e_j)\xi,f_i\big>=e_j(\lambda_j)\delta_{ij}
    +\lambda_jF_{j| ji}-\lambda_i E_{j| ji}, \\[1ex]
    \big<r(e_i,e_j)\xi,f_j\big>=e_i(\lambda_j), \\[1ex]
    \big<r(e_i,e_j)\xi,f_i\big>=e_i(\lambda_j)\delta_{ij}+\lambda_j F_{i| ji}-
    \lambda_i E_{i| ji}
    \end{array}
    $$
    From this it follows that
    $$
    \begin{array}{l}
    \big<R(e_i,e_j)\xi,f_j\big>=\big<r(e_i,e_j)\xi,f_j\big>-\big<r(e_j,e_i)\xi,f_j\big>=\\[1ex]
    e_i(\lambda_j)-e_i(\lambda_j)\delta_{ij}-\lambda_j F_{i| ji}+
    \lambda_i E_{i| ji} =\\[1ex]
     e_i(\lambda_j) - e_j(\lambda_j)\delta_{ij}-\lambda_jF_{j| ji}+\lambda_i E_{j| ji}+
    (\lambda_i+\lambda_j)(E_{j| ij}-F_{j| ij})= \\[1ex]
     e_i(\lambda_j)-\big<r(e_j,e_j)\xi,f_i\big>-
    (\lambda_i+\lambda_j)(E_{j| ji}-F_{j| ji}).
    \end{array}
    $$
    So, we see that
    $$
    \big<r(e_j,e_j)\xi,f_i\big>=e_i(\lambda_j)-
    (\lambda_i+\lambda_j)(E_{j| ji}-F_{j| ji})
    -\big<R(e_i,e_j)\xi,f_j\big>.
    $$
    Finally, introducing the matrix $G_{i| j}$ with the components
    $$
    G_{i| j}=E_{i| ij}-F_{i| ij},
    $$
    we can rewrite the mean curvature formula as follows

\begin{equation}\label{SH}
\begin{array}{c}
\displaystyle    (n+1)H_{\sigma |}=\\[2ex]
\displaystyle  \frac{1}{\sqrt{1+\lambda_\sigma^2}}
    \sum_{i=0}^n\frac{e_\sigma(\lambda_i)-(\lambda_i+\lambda_\sigma)G_{i| \sigma}
    +(\lambda_i\lambda_\sigma-1)\big<R(e_\sigma,e_i)\xi,f_i\big>}
    {1+\lambda_i^2},
\end{array}
\end{equation}

  where $\lambda_0=0$ and $f_0=0$ is supposed.

\section{Some special cases and examples}

\subsection{Normal vector field of a Riemannian foliation}

    We consider an important special case of a unit {\it geodesic } vector field
    $\xi$ such that the orthogonal distribution $\xi^\perp$ is integrable.
    In other words, suppose that a given Riemannian manifold admits a Riemannian
    transversally orientable hyperfoliation. Then the following holds.

\begin{theorem} Let $M^{n+1}$ admit a Riemannian transversally orientable
    hyperfoliation. Let $\xi$ be a unit normal vector field of the foliation. Then
    the components of the mean curvature vector of $\xi(M)$ are
    $$
    H_{\sigma |}=\frac{1}{(n+1)\sqrt{1+k_\sigma^2}}
    \sum_{\alpha=1}^{n}\left\{\frac{-e_\sigma(k_\alpha)+
    (1-k_\alpha k_\sigma)\big<R(\xi,e_\alpha)e_\alpha,e_\sigma\big>}{1+k_\alpha^2}
    \right\}
    $$
    where  $e_\alpha$ determine the principal directions and  $k_\alpha $ are the principal
    curvatures of the fibers.
\end{theorem}

    \begin{remark} \rm The analogous problem was treated in \cite{BX-V1}, where the
    authors considered the  {\it minimality} condition for the vector
    field. The corresponding conditions in \cite{BX-V1} differ from the mean curvature
    components by a factor. We refer to \cite{BX-V4} for applications of this conditions.
    \end{remark}

\begin{proof} For the given situation, the singular frame is simple. As $\xi$ is geodesic
    vector field, we have $e_0=\xi $, while the others are principal vectors of the second
    fundamental form of the fibers. If we denote the corresponding shape operator  by $A_\xi$,
    then
    $$
    \nabla_{e_\alpha}\xi=-A_\xi e_\alpha=-k_\alpha e_\alpha
    $$
    So, neglecting the condition on the $\lambda_\alpha $ to be {\it positive} (in fact, we never used
this condition
    in proof of the formula (\ref{SH})), we may put $f_\alpha=e_\alpha $ and $\lambda_\alpha=-
k_\alpha$.
    Therefore, in (\ref{SH}) we obtain $G_{i| j}=0$ and the result follows immediately.

\end{proof}

\subsection{Strongly normal vector field.}

A unit vector field $\xi$ is called  {\it normal} if
$R(X,Y)\xi=\alpha \xi $ and {\it strongly normal} if
$r(X,Y)\xi=\alpha \xi$ for all $X,Y \in \xi^\perp $. Our result
(\ref{H}) allows to prove easily \cite{GD-V1}:

\vspace{1ex}
 { \it Every unit strongly normal geodesic vector
field is minimal}

\vspace{1ex}

Indeed, since $\xi$ is geodesic, $\nabla_\xi\xi=0$ and therefore
$e_0=\xi$. Hence, $r(e_0,e_0)\xi=0$ and $e_1, \dots , e_{n} \in
\xi^\perp , \quad f_1, \dots , f_{n} \in \xi^\perp$. Evidently, a
strongly normal vector field is always normal. So, each term in
(\ref{H}) vanishes.

\subsection{Geodesic vector fields on 2-dimensional manifolds}

    For $dim M=2$ the mean curvature of $\xi(M) \subset T_1M$ equals
    $$
    H= \frac{1}{2 \sqrt{1+\lambda^2}} \left\{ \big< r(e_0,e_0) \xi +
    \frac{r(e_1,e_1) \xi}{1+\lambda^2},f_1 \big> \right\}
    $$
    or
    \begin{equation}\label{H_2}
        H=\frac{1}{2\sqrt{1+\lambda^2}}\left\{-\big<\nabla_{e_0}e_0,e_1\big>
        \lambda+\frac{e_1(\lambda)}{1+\lambda^2}\right\}.
    \end{equation}

    The above formula allows to prove the following statement.
    \vspace{1ex}

    {\it A unit geodesic vector field on a 2-dimensional manifold is
    minimal if and only if it is strongly normal} (see \cite{GD-V1}).
    \vspace{1ex}

    Indeed, in this case we can set $e_0=\xi$, $f_1=\pm e_1$. So,
    up to a sign,
    $$
     H=\frac{1}{2(1+\lambda^2)^{3/2}}\big<r(e_1,e_1)\xi,e_1\big>
    $$
     and the statement follows immediately.

    In \cite{GD-V1}, the authors give an example of a geodesic but not
    strongly normal vector field and hence not minimal. Here we can
    easily find the mean curvature of that field. Namely, consider the
    2-dimensional manifold of non-positive curvature with metric
    $$
    ds^2=du^2+e^{2uv}dv^2.
    $$
    Set $\xi=\{1,0\}$. Then, up to a sign, the
    singular frame is
    $$
    e_0=\xi \mbox{  and  } e_1=\{0,e^{-uv}\}=f_1.
    $$
    It is easy to see that $$ \nabla _{e_1}\xi=ve_1. $$ Hence
    $\lambda=v$ and  $e_1(\lambda)=e^{-uv}$.
    So, the mean curvature of $\xi(M)$ is given by
    $$
    H=\frac{e^{-uv}}{2(1+v^2)^{3/2}}.
    $$

    \vspace{1ex}

\subsection{Examples of non-geodesic minimal vector fields on some 2-dimensional
            Riemannian manifolds}

    Next, we consider a Riemannian 2-manifold $M$ with the metric
    $$
    ds^2=du^2+e^{2g(u)}dv^2.
    $$
    As it was shown in \cite{GD-V1} for
    the general situation, the vector field $\partial/\partial u$ is
    minimal. Here we shall consider the vector field which makes a
    constant angle with $\partial/\partial u$ along  each $u$ - geodesic.

    \begin{proposition}Up to a sign,
    the mean curvature of the vector field $\xi$ on a 2-dimensional
    Riemannian manifold with metric $ ds^2=du^2+e^{2g(u)}dv^2$ which
    is parallel along each $u$ - geodesic, is
    $$
    H=\frac{e^{-2g}\omega_{vv}}{2\Big(1+(e^{-g}\omega_v+g')^2\Big)^{3/2}},
    $$
    where $\omega(v)$ is the angle function of $\xi$ with respect
    to the direction of $u$ - geodesics.
    \end{proposition}

    \begin{proof}
    Consider the mutually orthogonal unit vector fields $ X_1= \{ 1,0
    \} $ and  $ X_2= \{ 0,e^{-g} \} $. A direct calculation gives $$
    \begin{array}{lclclcl}
    \nabla_{X_1}X_1 & = & 0,          &  & \nabla_{X_1}X_2 & = & 0 ,\\
    \nabla_{X_2}X_1 & = & g' X_2, &  & \nabla_{X_2}X_2 & = & -g'X_1.
    \end{array}
    $$ Let $\omega (u,v)$ be the angle function defining the vector
    field $\xi$ by $$ \xi=\cos \omega X_1 + \sin \omega X_2 $$ Let
    $\eta $ be a unit vector field orthogonal to $\xi$: $$ \eta=-\sin
    \omega X_1 + \cos \omega  X_2. $$ Then $$ \nabla_{X_1}\xi =
    X_1(\omega) \eta , \mbox{       }
    \nabla_{X_2}\xi=-(X_2(\omega)+g')\eta. $$

    Now, suppose $\xi$ to be parallel along a $u$ - geodesic, that is,
    set $X_1(\omega)=0$. Then the singular frame is : $e_0=X_1$ and $
    e_1=X_2$. The singular function is $\lambda=-(X_2(\omega)+g')$ and
    we see that, up to a sign, $f_1$ coincides with $\eta$. So
    $$
    H=\frac{e_1(\lambda )}{2(1+\lambda^2)^{3/2}}.
    $$
    For $e_1(\lambda)$ we obtain
    $$
    e_1(\lambda)=X_2(-X_2(\omega)+g')=-X_2(X_2(\omega))+X_2(g')=-e^{-2g}\omega_{vv}
    $$
     since $g$ does not depend on $v$. Therefore
    $$
    H=\frac{e^{-2g}\omega_{vv}}{2\Big(1+(e^{-g}\omega_v+g')^2\Big)^{3/2}},
    $$
    what was claimed.
    \end{proof}

    From the above formula we conclude:

    \vspace{1ex}
    {\it On a 2-dimensional manifold with metric $ds^2=du^2+e^{2g(u)}dv^2$ the unit
    vector field  $\xi$which is parallel along $u$ -- geodesics, is minimal if its angle
    increment along $v$ -- curves is not higher then the linear one. }
    \vspace{1ex}

    Particularly, if $\omega=const$, then $\xi$ is minimal.

\subsection{The mean curvature of a general unit vector field on
            2-dimensional manifolds}
    In the case of $dim M=2,$ the mean curvature of a unit vector field can
    be expressed in terms of the geodesic curvature of integral curves
    of the given field and their orthogonal trajectories.
\begin{proposition}
    Let $\xi$ and $\eta$ be unit mutually orthogonal vector fields
    on a 2-dimensional Riemannian manifold. Denote by $k$ and $\kappa$
    the geodesic curvatures of the integral curves of the field $\xi$
    and $\eta$, respectively. The mean curvature $H$ of the vector
    field $\xi$ is given, up to a sign, by
    $$
      H=\frac12\left[\xi\left(\frac{k}{\sqrt{1+k^2+\kappa^2}}\right)-
      \eta\left(\frac{\kappa}{\sqrt{1+k^2+\kappa^2}}\right)
      \right].
    $$

\end{proposition}

\begin{remark} \rm  The analogous expression can be found in
                \cite{GM-LF} as a condition of minimality of the unit vector field
                on 2-dimensional manifolds.
\end{remark}

\begin{proof}
    From (\ref{H_2}) one can see that after the replacement $\xi\to -\xi$
    the mean curvature $H$ just changes its sign. Therefore, we may
    choose the direction of $\xi$ in such a way that it will be the
    field of principal normals of the $\eta $ -- curves. The same arguments
    allow us to consider $\eta$ as the field of principal normals
    of the $\xi$ -- curves. Denote by $\omega$ an angle between $\xi$ and the
    field $e_0$ of the singular frame. Then

    $$ e_0=\cos\omega\xi+ \sin\omega\eta. $$
    As $\nabla_{e_0}\xi=0$, we have
    $$\cos\omega\nabla_\xi\xi+\sin\omega\nabla\eta\xi=0.$$
    The Frenet formulas give
    $$\nabla_\xi\xi=k\eta, \quad \nabla_\eta\xi=-\kappa\eta.$$
    Therefore, we obtain
\begin{equation}
\label{par}
            k\cos\omega-\kappa\sin\omega=0.
\end{equation}
        Denote by $e_1$ and $f_1$ the other vectors of the singular
        frame. It is easy to check that the change of directions of these
        vectors induces a sign change of $H$. Therefore, we can always
        set $f_1=\eta$ and $e_1=\pm\sin\omega\xi \mp \cos\omega\eta$
        to satisfy the equation $\nabla_{e_1}\xi=\lambda f_1$ with
        $\lambda\geq 0$. Taking all of this into account, set
 $$
 \begin{array}{c}
     e_0=\cos\omega\xi+\sin\omega\eta, \\
     e_1=\sin\omega\xi-\cos\omega\eta.
 \end{array}
 $$
       Then we have
$$
\begin{array}{l}
        \nabla_{e_0}\xi=\cos\omega\nabla_\xi\xi
        +\sin\omega\nabla_\eta\xi=0, \\[1ex]
         \nabla_{e_1}\xi=\sin\omega\nabla_\xi\xi-\cos\omega\nabla_\eta\xi=\lambda\eta.
\end{array}
$$
    From these equations we derive
$$
\begin{array}{l}
        \nabla_\xi\xi=\lambda\sin\omega \,\eta,\\[1ex]
        \nabla_\eta\xi=-\lambda\cos\omega\,\eta.
\end{array}
$$
    Comparing this with the Frenet formulas, we conclude that
    $k=\lambda\sin\omega,\ \kappa=\lambda\cos\omega$. Therefore,
\begin{equation}
\label{om}
        \lambda^2=k^2+\kappa^2,\quad \sin\omega=\frac{k}{\lambda},
        \quad \cos\omega=\frac{\kappa}{\lambda}
\end{equation}

    To use the formula (\ref{H_2}), we should find $e_1(\lambda)$ and
    $\big<\nabla_{e_0}e_0,e_1\big>$. Now, keeping in mind (\ref{par}), we have

    $$  e_1(\lambda)
        =\frac{k}{\lambda}\xi(\lambda)-\frac{\kappa}{\lambda}\eta(\lambda)
    $$
    and
$$
\begin{array}{rl}\displaystyle
    \nabla_{e_0}e_0=&\cos{\omega}\nabla_\xi(\cos{\omega}\,\xi+\sin{\omega}\,\eta)
    +\sin{\omega}\nabla_\eta(\cos{\omega}\,\xi
    +\sin{\omega}\,\eta)=\\\displaystyle
    &-(\xi(\omega)\cos\omega+\eta(\omega)\sin\omega)e_1-
    (k\cos\omega-\kappa\sin\omega)e_1=\\\displaystyle
    &-\big(\xi(\sin\omega)-\eta(\cos\omega)\big)e_1.
\end{array}
$$
    Therefore, using (\ref{om}), we get
    $$
        -\big<\nabla_{e_0}e_0,e_1\big>=\xi\left(\frac{k}{\lambda}\right)-
        \eta\left(\frac{\kappa}{\lambda}\right).
    $$
    Substituting these  expressions into (\ref{H_2}), we
    obtain
$$
\begin{array}{l}
    \displaystyle H=\!\frac12\frac{1}{\sqrt{1+\lambda^2}}
    \left[\left(\xi\Big(\frac{k}{\lambda}\Big)-
    \eta\Big(\frac{\kappa}{\lambda}\Big)\right)\lambda+
    \frac{1}{1+\lambda^2}\left(\frac{k}{\lambda}\xi(\lambda)-
    \frac{\kappa}{\lambda}\eta(\lambda)\right)\right]=\\[2ex]\displaystyle
    \qquad\frac{1}{2}\frac{1}{(1+\lambda^2)^{3/2}}
    \left[\big((1+\lambda^2)\,\xi(k)-k\lambda\,\xi(\lambda)\big)-
    \big((1+\lambda^2)\,\eta(\kappa)-\kappa\lambda\,\eta(\lambda)\big)\right]=\\[3ex]\displaystyle
    \qquad\frac12\left[\xi\left(\frac{k}{\sqrt{1+\lambda^2}}\right)-
    \eta\left(\frac{\kappa}{\sqrt{1+\lambda^2}}\right)\right].
\end{array}
$$
    Taking into account (\ref{om}), we get what was claimed.
\end{proof}
\vspace{1ex}

    {\bf Corollary. }{\it If $\xi$ is a geodesic vector field  then

    $$
    H=-\frac12\frac{\partial}{\partial\sigma}\left(\frac{\kappa}{\sqrt{1+\kappa^2}}\right)
    $$
    where $\sigma$ is the arc-length parameter of the orthogonal trajectories of the field $\xi$
    and $\kappa$ is their geodesic curvature.}
    \vspace{1ex}

    A unit geodesic vector field is said to be {\it radial} if it is a tangent vector field of
    geodesics starting at a fixed point. Now we can confirm the following statement \cite{BX-V1}.
\begin{proposition}
        If each radial vector field on a 2-dimensional Riemannian
        manifold $M$ is minimal, then $M$ has constant curvature.
\end{proposition}

\begin{proof}
        Indeed, if such a vector field is minimal, then its orthogonal
        trajectories are Gauss circles of constant geodesic curvature,
        which means that those circles are Darboux ones. Therefore,
        $M$ is of constant Gaussian curvature ( see \cite{Bl}).
\end{proof}

\subsection{Some examples of vector fields of constant mean curvature.}
\subsubsection{The example on the Lobachevsky 2-space.}
    Consider the Lobachevsky plane $L^2$ with the metric
    $$
    ds^2=du^2+e^{2u}dv^2.
    $$
    The coordinate lines of $L^2$ are $u$ -geodesics and their orthogonal trajectories.

\begin{proposition}
\label{Ex}
    The unit vector field on $L^2$ whose angle function with respect to
    $u$ - geodesics is $\omega=au+b \ (a,b=const)$ has constant mean
    curvature
    $$
    H=\frac{a}{2\sqrt{2+a^2}}.
    $$
\end{proposition}

\begin{proof}
    Indeed, consider the field $\xi=\cos\omega X_1+\sin\omega X_2$ where
    $\omega=au+b$ and $X_1=\{1,0\}, \ X_2=\{0,e^{-u}\}$. Then
$$
\begin{array}{ll}
    \nabla_{X_1}X_1=0,    & \nabla_{X_1}X_2=0, \\[1ex]
    \nabla_{X_2}X_1=X_2,  & \nabla_{X_2}X_2=-X_1.
\end{array}
$$
    Now we define the singular frame for $\xi$. To do this, we
    introduce the vector field  $\eta=-\sin\omega X_1+\cos\omega X_2$.
    Then
$$
\begin{array}{l}
    \nabla_{X_1}\xi=\frac{\partial\omega}{\partial u}\eta=a\eta ,\\[1ex]
    \nabla_{X_2}\xi=\eta.
\end{array}
$$
    Therefore, setting
    $$
    e_0=\frac{1}{\sqrt{1+a^2}} ( X_1-aX_2), \ \ \
    e_1=\frac{1}{\sqrt{1+a^2}} (a X_1+X_2),
    $$
    we have
    $$
    \nabla_{e_0}\xi=0, \ \ \ \nabla_{e_1}\xi=\sqrt{1+a^2}\eta.
    $$
    Hence, $f_1=\eta$ and $ \lambda=\sqrt{1+a^2}=const$. So,
    $e_1(\lambda)=0$. Moreover,
    $$
    \nabla_{e_0}e_0=-\frac{a}{\sqrt{1+a^2}}\,e_1.
    $$
    Substituting this into  (\ref{H_2}), we have
    $$
    H=\frac{a}{2\sqrt{2+a^2}}.
    $$
    So, the statement is proved.
\end{proof}

\subsubsection{The generalized examples on the Lobachevsky  $(n+1)$- space.}

    Consider the $(n+1)$ - dimensional Lobachevsky space endowed with
    horospherical coordinates $( u, v^1, \dots, v^n)$. Then
    $$
     ds^2=du^2+e^{2u} [(dv^1)^2+ \dots + (dv^n)^2 ].
    $$

    Consider the unit vector fields
    \begin{equation}\label{base}
        X_0=\{1,0,\dots, 0\}, X_1=\{0,e^{-u},\dots,0\},\dots, X_n=\{0,0,\dots,e^{-u}\}.
    \end{equation}
    It is easy to check that
    $$
    \begin{array}{ll}
    \nabla_{X_{\scriptstyle 0}}X_0=0, &\nabla_{X_{\scriptstyle 0}}X_{\alpha}=0,\\
    \nabla_{X_{\scriptstyle \alpha}}X_0=X_{\alpha} &
    \nabla_{X_{\scriptstyle\alpha}}X_\alpha=-X_0.
    \end{array}
    $$
    Define the unit vector field $\xi$ as follows:
\begin{equation}\label{vf1}
    \xi=\cos\theta X_0+\sin\theta\cos u X_1+\sin\theta\sin uX_2,
\end{equation}
    where $\theta\in [0, \pi/2]$ is constant.

\begin{proposition}\label{Lob_n}
    The unit vector field which is given by (\ref{vf1}) with respect to the
    frame (\ref{base}) on Lobachevsky $(n+1)$ - space with the metric
    $$
     ds^2=du^2+e^{2u} [(dv^1)^2+ \dots + (dv^n)^2 ],
    $$
    is a field of constant mean curvature. Namely, we have
    $$
    \begin{array}{l}
    H_{1|}=\displaystyle\frac{n-2}{n+1}\frac{\sqrt{2}\sin\theta\cos\theta}{1+\cos^2{\theta}},\\[2ex]
    H_{2|}=\displaystyle\frac{n\sqrt{2}\sin\theta}{2(n+1)}, \\[2ex]
    H_{\sigma|}=0 \quad \sigma\geq 3.
    \end{array}
    $$
\end{proposition}

\begin{proof}

    With respect to the frame $\{X_0,X_1,\dots,X_n\}$, the matrix
    $(\nabla\xi)$ has the form
    $$
    \left[
    \begin{array}{cccccc}
    0                    &   -\sin\theta\cos u  &  -\sin\theta \sin u &0&\dots&0\\
    -\sin\theta\sin u &\cos\theta                &        0               &0&\dots&0\\
    \sin\theta\cos u  & 0                       &  \cos\theta             &0&\dots&0\\
         0                   & 0                       &    0                   &\cos\theta    & \dots&0 \\
    \vdots                &\vdots                  &  \vdots                &0&\ddots&0\\
    0                    &0                        & 0                      &0&\dots &\cos\theta
        \end{array}
    \right].
    $$
    It is easy to find that the matrix $(\nabla\xi)^t(\nabla\xi)$ has the following
    expression
    $$
    \left[
    \begin{array}{cc}
        A&0\\
        0&B
    \end{array}
    \right],
    $$
    where $A$ is the $3\times 3$ matrix
    $$
    \left[
    \begin{array}{ccc}
    \sin^2\theta                    &   -\sin\theta\cos\theta\sin u  &  \sin\theta\cos\theta\cos u\\
    -\sin\theta\cos\theta\sin u &\cos^2\theta +\sin^2\theta\cos^2u
&\sin^2\theta\sin u\cos u \\
    \sin\theta\cos\theta\cos u  & sin^2\theta\sin u\cos u & \cos^2\theta
   +\sin^2\theta\sin^2(u)\\
        \end{array}
    \right]
    $$
    and $B$ is the diagonal $(n-2)\times(n-2)$ matrix of the form
    $$
    \left[
    \begin{array}{ccc}
        \cos^2\theta &\dots&0\\
        \vdots        &\ddots&\vdots \\
        0             & \dots &\cos^2\theta
    \end{array}
    \right].
    $$
    The eigenvalues of the matrix $(\nabla\xi)^t(\nabla\xi)$ are
    $$
    \lambda_0^2=0, \lambda_1^2=\lambda_2^2=1,
    \lambda_3=\dots=\lambda_n^2=cos^2\theta.
    $$
    Now it is easy to find the vectors of the singular frame. We get
    $$
    \begin{array}{l}
    \begin{array}{ccl}
    e_0&=&\cos\theta  X_0+\sin\theta\sin uX_1-\sin\theta\cos uX_2, \\
    e_1&=&\cos uX_1+\sin uX_2, \\
    e_2&=&\sin\theta X_0-\cos\theta  \sin uX_1+\cos\theta  \cos uX_2, \\
    \end{array}\\
    \ \, e_3= X_3,\dots , e_n= X_n
    \end{array}
    $$
    and
    $$
    \begin{array}{l}
    \begin{array}{ccl}
    f_1&=&-\sin\theta X_0+\cos\theta  \cos u X_1+\cos\theta\sin uX_2,\\
    f_2&=&-\sin u X_1+\cos uX_2, \\
    \end{array}\\
    \ \, f_3=e_3, \dots , f_n =e_n.
    \end{array}
    $$
    So, we have
    $$
    \begin{array}{lcl}
    \nabla_{\displaystyle e_0}\xi=0, & \nabla_{\displaystyle e_1}\xi=f_1,
    & \nabla_{\displaystyle e_2}\xi= f_2,\\[1ex]
    \nabla_{\displaystyle e_3}\xi=\cos\theta  f_3, & \dots
    & \nabla_{\displaystyle e_n}\xi=\cos\theta f_n
    \end{array}
    $$

    Straightforward computation gives the following components for the
    matrix $G_{i| j}$:
    $$
    \left[
    \begin{array}{cccccc}
    0       & \sin\theta\cos\theta   & -\sin\theta & 0 & \dots & 0 \\
    -\cos\theta   & 0                  & -\sin\theta & 0 & \dots & 0 \\
    -\cos\theta   & -\sin\theta\cos\theta   & 0 & 0 &\dots & 0 \\
       -\cos\theta   & -\sin\theta & -\sin\theta &0 & \dots & 0 \\
    \vdots & \vdots & \vdots & \vdots & \vdots & \vdots \\
    -\cos\theta   & -\sin\theta & -\sin\theta &0 & \dots & 0
    \end{array}
    \right].
    $$
    As all $ \lambda_i$ are constants, we have
    $$
    \begin{array}{cl}
    H_{1|}=&\displaystyle\frac{1}{(n+1)\sqrt{1+\lambda_1^2}}
    \sum\limits_{i=0}^{n}\frac{-(\lambda_1+\lambda_i)G_{i| 1} +
    (\lambda_i-\lambda_1)\big<R(e_1,e_i)\xi,f_i\big>}{1+\lambda_i^2}=\\[3ex]
    &\displaystyle \frac{1}{(n+1)\sqrt{2}}\left[
    \sum\limits_{i=0}^2 (-G_{i| 1})+\sum\limits_{i=3}^n
    \frac{-(1+\lambda_i)G_{i| 1}+(\lambda_i-1)\big<\xi,e_1\big>}
        {1+\cos^2\theta}\right] = \\[3ex]
    &\displaystyle\frac{1}{(n+1)\sqrt{2}}\left[0+(n-2)\frac{(1+\cos\theta  \sin\theta+(\cos\theta  -
    1)\sin\theta} {1+\cos^2\theta}\right]= \\[3ex]
    &\displaystyle\frac{n-2}{n+1}\frac{\sqrt{2}\sin\theta\cos\theta  }{1+\cos^2\theta }.
    \end{array}
    $$
    Analogously, we get
    $$
    \begin{array}{cl}
    H_{2|}=&\displaystyle\frac{1}{(n+1)\sqrt{1+\lambda_2^2}}
    \sum\limits_{i=0}^{n}\frac{-(\lambda_2+\lambda_i)G_{i| 2} +
    (\lambda_i-\lambda_2)\big<R(e_2,e_i)\xi,f_i\big>}{1+\lambda_2^2}=\\[3ex]
     &\displaystyle\frac{1}{(n+1)\sqrt{2}}\left[
    \sum\limits_{i=0}^2 (-G_{i| 2})+\sum\limits_{i=3}^n
    \frac{-(1+\lambda_i)G_{i| 2}+(\lambda_i-1)\big<\xi,e_2\big>}
        {1+\cos^2\theta}\right] = \\[3ex]
     &\displaystyle\frac{\sqrt{2}}{2(n+1)}\left[2\sin\theta+(n-2)
    \frac{(1+\cos\theta  )\sin\theta+(\cos\theta  -1)\sin\theta\cos\theta  }
        {1+\cos^2\theta}\right]=\\[3ex]
     &\displaystyle\frac{\sqrt{2}}{2(n+1)}\left[2\sin\theta+(n-2)
    \frac{\sin\theta+\sin\theta\cos^2\theta}
        {1+\cos^2\theta}\right]=
    \frac{n\sqrt{2}\sin\theta}{2(n+1)}.
    \end{array}
    $$
    and $H_{\sigma|}=0$ for all $\sigma\geq 3$.
\end{proof}

    A similar but more complicated computation shows that there exist a
    family of vector fields of constant mean curvature on the Lobachevsky space.
    Namely, let $\xi$ be a vector field given by
\begin{equation} \label{vf2}
    \xi=\cos\theta   X_0+\sin\theta\cos{au} X_1+\sin\theta\sin{au} X_2,
\end{equation}
    where $a$ and $\theta$ are constants and the frame $X_0, X_1, \dots , X_n$
    is chosen as above. Then the following statement is true.

\begin{proposition}\label{Lob_n2}
    The unit vector field which is given by (\ref{vf2}) with respect to the
    frame (\ref{base}) on the Lobachevsky $(n+1)$ - space with the metric
    $$
     ds^2=du^2+e^{2u} [(dv^1)^2+ \dots + (dv^n)^2 ],
    $$
    is a field of constant mean curvature. Namely, we have
    $$
    \begin{array}{l}
    H_{1|}=\displaystyle\frac{\sqrt{2}\sin\theta\cos\theta  }{n+1}
    \left(\frac{1-a^2}{1+\cos^2\theta+a^2\sin^2\theta}+\frac{n-
    2}{1+\cos^2\theta}\right),\\[2ex]
    H_{2|}=\displaystyle\frac{an\sin\theta}{(n+1)\sqrt{1+\cos^2\theta+a^2\sin^2\theta}},\\[2ex]
    H_{\sigma|}=0 \quad \sigma\geq 3.
    \end{array}
    $$
\end{proposition}

    The proof is based on the fact that the singular values of $(\nabla\xi)$
    are the following constants:
    $$
    \lambda_1=1, \lambda_2=\sqrt{\cos^2\theta+a^2\sin^2\theta}, \lambda_3=
    \ldots =\lambda_n=\cos\theta  .
    $$
    {\large Acknowledgement.} The author expresses his thanks to
    P.Nagy who invited him to take part in a fruitful workshop
    (Debrecen, 2000) on the geometry of tangent sphere bundle. The talks with L.Vanhecke and
    E.Boeckx gave the starting impulse to the article. The author
    also thanks A.Borisenko who was the first who asked on
    examples of vector fields of constant mean curvature.

\end{document}